\documentclass{pnastwo_custom}

\usepackage{amsmath}
\usepackage{graphicx,epsf}
\usepackage{amssymb,amsfonts}

\newcommand{\eg}{{\em e.g.}}

\newcommand{\cf}{{\em cf. }}

\newcommand{\real}{{\mathbb R}}

\newcommand{\nats}{{\mathbb N}}

\newcommand{\zed}{{\mathbb Z}}

\newcommand{\abs}[1]{\left\vert{#1}\right\vert}

\newcommand{\inv}{^{-1}}
\newcommand{\one}{{\sf{1}}}    
\newcommand{\Support}{{\mathcal S}}     
\newcommand{\Index}{{\mathcal I}}       
\newcommand{\style}[1]{{\sc{#1}}}    
\newcommand{\Nodes}{{\mathcal N}}       
\newcommand{\Transform}{{\mathcal T}}         
\newcommand{\Triangulation}{{\mathcal T}}         
\newcommand{\crit}{{\mathcal C}}        
\newcommand{\Def}{{\sf Def}}            
\newcommand{\id}{{\sf Id}}              
\newcommand{\totvar}{{\sf totvar}}      
\newcommand{\Dual}{{\mathcal D}}
\newcommand{\dchifloor}{{\lfloor d\chi\rfloor}}
\newcommand{\dchiceil}{{\lceil d\chi\rceil}}
\newcommand{\qed}{{$\blacksquare$}}
\newcommand{\Omin}{{\mathcal O}}

\begin{document}
%

\title{Euler integration over definable functions}

\author{
Yuliy Baryshnikov\affil{1}{Bell Laboratories, Murray Hill, NJ},
Robert Ghrist\affil{2}{University of Pennsylvania, Philadelphia, PA}
}

\contributor{
{\bf Preprint version: August 2009.}

}

\maketitle


\begin{article}

\begin{abstract}
We extend the theory of Euler integration from the class of
constructible functions to that of ``tame'' $\real$-valued functions
(definable with respect to an o-minimal structure). The corresponding
integral operator has some unusual defects (it is not a linear
operator); however, it has a compelling Morse-theoretic
interpretation. In addition, we show that it is an appropriate setting
in which to do numerical analysis of Euler integrals, with applications
to incomplete and uncertain data in sensor networks.
\end{abstract}


\dropcap{I}ntegration with respect to Euler characteristic is a
homomorphism $\int_X\cdot d\chi:CF(X)\to\zed$ from the ring of
constructible functions $CF(X)$ (``tame'' integer-valued functions on
a topological space $X$) to the integers $\zed$. It is a topological
integration theory which uses as a measure the venerable Euler
characteristic $\chi$. Euler characteristic, suitably defined, satisfies
the fundamental property of a measure:
\begin{equation}
\label{eq:measure}
    \chi(A\cup B) = \chi(A) + \chi(B) - \chi(A\cap B),
\end{equation}
for $A$ and $B$ ``tame'' subsets of $X$. We extend the theory to
$\real$-valued integrands and demonstrate its utility in managing
incomplete data in, e.g., sensor networks.

\section{Constructible integrands}
\label{sec:z-valued}

Because the Euler characteristic is only finitely additive, one must
continually invoke the word ``tame'' to ensure that $\chi$ is
well-defined. One means by which to do so it via an \style{o-minimal
structure} \cite{vdD}, a sequence $\Omin=(\Omin_n)$ of Boolean
algebras of subsets of $\real^n$ satisfying a small list of axioms:
closure under products, closure under projections, and finite
decompositions in $\Omin_1$. Elements of $\Omin$ are called
\style{definable} sets and these are ``tame'' for purposes of
integration theory. Examples of o-minimal structures include (1)
piecewise-linear sets;\footnote{Some authors require an o-minimal
structure to contain algebraic curves, eliminating this particular
example.} (2) semi-algebraic sets; and (3) globally subanalytic sets.

Definable functions between spaces are those whose graphs are in
$\Omin$. For $X$ and $Y$ definable spaces, let $\Def(X,Y)$ denote
the class of compactly supported definable functions $h:X\to Y$, and
fix as a convention $\Def(X)=\Def(X,\real)$. Let
$CF(X)=\Def(X,\zed)\subset\Def(X,\real)$ denote the ring of
\style{constructible functions}: compactly supported $\zed$-valued
functions all of whose level sets are definable. Note that in general,
definable functions (even definable `homeomorphisms') are not
necessarily continuous.


We briefly recall the theory of Euler integration, established as an
integration theory in the constructible setting in
\cite{MacPherson,Schapira,Schapira:tom,Viro} and anticipated by a
combinatorial version in \cite{Blaschke,Groemer,Hadwiger,Rota}. Fix
an o-minimal structure $\Omin$ on a space $X$. The geometric Euler
characteristic is the function $\chi:\Omin\to\zed$ which takes a
definable set $A\in\Omin$ to $\chi(A)=\sum_i(-1)^i\dim
H_i^{BM}(A;\real)$, where $H_*^{BM}$ is the Borel-Moore homology
(equivalently, singular compactly supported cohomology) of $A$. This
also has a combinatorial definition: if $A$ is definably homeomorphic
to a finite disjoint union of (open) simplices $\coprod_j \sigma_j$,
then $\chi(A)=\sum_j(-1)^{\dim\sigma_j}$. Algebraic topology asserts
that $\chi$ is independent of the decomposition into simplices. The
Mayer-Vietoris principle asserts that $\chi$ is a measure (or
`valuation') on $\Omin$, as expressed in Eqn. [\ref{eq:measure}].

The \style{Euler integral} is the pushforward of the
trivial map $X\mapsto \{pt\}$ to $\int_Xd\chi:CF(X)\to
CF(\{pt\})\cong\zed$ satisfying $\int_X\one_A\,d\chi=\chi(A)$ for
$\one_A$ the characteristic function over a definable set $A$. From
the definitions and a telescoping sum one easily obtains:
\begin{equation}
\label{eq:level}
\int_X h\,d\chi
    =
    \sum_{s=-\infty}^\infty s \chi\{h=s\}
    =
    \sum_{s=0}^\infty \chi\{h>s\} - \chi\{h<-s\}
.
\end{equation}
Because the Euler integral is a pushforward, any definable map
$F:X\to Y$ induces $F_*:CF(X)\to CF(Y)$ satisfying $\int_Xh\,d\chi =
\int_YF_*h\,d\chi$. Explicitly,
\begin{equation}
\label{eq:fubini}
    F_*h(y)
    =
    \int_{F\inv(y)}h\,d\chi,
\end{equation}
as a manifestation of the Fubini Theorem.

The Euler integral has been found to be an elegant and useful tool for
explaining properties of algebraic curves \cite{Brocker} and
stratified Morse theory \cite{Schurmann,BK}, for reconstructing
objects in integral geometry \cite{Schapira:tom}, for target counting
in sensor networks \cite{BG:enum}, and as an intuitive basis for the more
general theory of motivic integration \cite{CL,Cluckers}.

\section{Real-valued integrands}
\label{sec:r-valued}

We extend the definition of Euler integration to $\real$-valued
integrands in $\Def(X)$ via step-function approximations.

\subsection{A Riemann-sum definition}

\begin{definition}
\label{def:realintegral}
Given $h\in \Def(X)$, define: 
\begin{eqnarray}
    \int_X h\,\dchifloor
    &=&
    \lim_{n\to\infty}
    \frac{1}{n}\int_X \lfloor nh\rfloor d\chi .
    \\
    \int_X h\,\dchiceil
    &=&
    \lim_{n\to\infty}
    \frac{1}{n}\int_X \lceil nh\rceil d\chi .
\end{eqnarray}
\end{definition}

We establish that these limits exist and are well-defined, though
not equal.

\begin{lemma}
\label{lem:r-simplex}
Given an affine function $h\in\Def(\sigma)$ on an open $k$-simplex
$\sigma$,
\begin{equation}
    \int_\sigma h\,\dchifloor = (-1)^k\inf_\sigma h
    \, ; \,
    \int_\sigma h\,\dchiceil = (-1)^k\sup_\sigma h    .
\end{equation}
\end{lemma}
{\em Proof:} For $h$ affine on $\sigma$, $\chi\{\lfloor nh\rfloor >
s\}=(-1)^k$ for all $s<n\inf_\sigma h$, and $0$ otherwise. One
computes
\[
    \lim_{n\to\infty}
    \frac{1}{n}\int_{\sigma}\lfloor nh\rfloor d\chi
    =
    \lim_{n\to\infty}
    \frac{1}{n}\sum_{s=0}^\infty \chi\{\lfloor nh\rfloor > s\}
    =
    (-1)^k\inf_\sigma h
    .
\]
The analogous computation holds with $\chi\{\lceil nh\rceil >
s\}=(-1)^k$ for all $s<n\sup_\sigma h$, and $0$ otherwise. \qed

This integration theory is robust to changes in coordinates.

\begin{lemma}
\label{lem:cov}
Integration on $\Def(X)$ with respect to $\dchifloor$ and
$\dchiceil$ is invariant under the right action of definable
bijections of $X$.
\end{lemma}
{\em Proof:} This is true for Euler integration on $CF(X)$; thus, it
holds for $\int_X\lfloor nh\rfloor\,d\chi$ and $\int_X\lceil
nh\rceil\,d\chi$. \qed

\begin{lemma}
\label{lem:welldef}
The limits in Definition \ref{def:realintegral} are well-defined.
\end{lemma}
{\em Proof:}
The \style{triangulation theorem} for $\Def(X)$ \cite{vdD} states that
to any $h\in\Def(X)$, there is a definable triangulation (a definable
bijection to a disjoint union of open affine simplices in some Euclidean
space) on which $h$ is affine
on each open simplex. The result now follows from Lemmas
\ref{lem:r-simplex} and \ref{lem:cov}. \qed

Integrals with respect to $\dchifloor$ and $\dchiceil$ are related to
total variation (in the case of compactly supported continuous
functions).

\begin{corollary}
\label{cor:totvar}
If $M$ is a 1-dimensional manifold and $h\in\Def(M)$ is continuous,
then
\begin{equation}
    \int_M h\,\dchifloor
    =
    -\int_M h\,\dchiceil
    =
    \frac{1}{2}\totvar(h) .
\end{equation}
\end{corollary}
{\em Proof:} Apply Lemma \ref{lem:r-simplex} to an affine
triangulation of $h$ which triangulates $M$ with the maxima
$\{p_i\}$ and minima $\{q_j\}$ as 0-simplices and the intervals
between them as 1-simplices. To each minimum $q_j$ is associated two
open 1-simplicies, since $M$ is a 1-manifold. Thus:
\[
    \int_M h\,\dchifloor = \sum_ih(p_i)+\sum_jh(q_j)-2\sum_jh(q_j)
    =\frac{1}{2}\totvar(h) .
\]
This equals $-\int_M h\,\dchiceil$ via an analogous computation.
\qed

This result generalizes greatly via Morse theory: see Corollary
\ref{cor:morse}. One notes that $\dchifloor$ and $\dchiceil$ give
integrals which are conjugate in the following sense.

\begin{lemma}
\label{lem:floortoceiling}
\begin{equation}
\label{eq:floortoceiling}
    \int h\dchiceil = -\int -h\dchifloor .
\end{equation}
\end{lemma}
{\em Proof:} Apply Lemma \ref{lem:r-simplex} to an affine
triangulation of $h$, and note that $\sup_\sigma h=-\inf_\sigma -h$.
\qed

The temptation to cancel the negatives must be resisted: see Lemma
\ref{lem:nonlinear} below.


\subsection{Computation}

Definition \ref{def:realintegral} has a Riemann-sum flavor which
extends to certain computational formulae. The following is a
definable analogue of Eqn. [\ref{eq:level}].

\begin{proposition}
\label{prop:r-valued}
For $h\in \Def(X)$,
\begin{eqnarray}
\label{eq:r-valued}
    \int_X h\,\dchifloor
    &=&
    \int_{s=0}^\infty \chi\{h\geq s\} - \chi\{h<-s\}\,ds
    \\
    \int_X h\,\dchiceil
    &=&
    \int_{s=0}^\infty \chi\{h> s\} - \chi\{h\leq -s\}\,ds
    .
\end{eqnarray}
\end{proposition}
{\em Proof:} For $h\geq 0$ affine on an open $k$-simplex $\sigma$,
\[
    \int_{\sigma} h\,\dchifloor
    = (-1)^k\inf_{\sigma} h
    = \int_0^\infty \chi(\sigma\cap\{h\geq s\})ds ,
\]
and for $h\leq 0$, the equation holds with $-\chi(\sigma\cap\{h <
-s\})$. The result for $\int\dchiceil$ follows from Lemma
\ref{lem:floortoceiling}. \qed

It is not true that $\int_Xh\,\dchifloor=\int_0^\infty s\chi\{h=s\}ds$:
the proper Lebesgue generalization of Eqn. [\ref{eq:level}] is the
following:

\begin{proposition}
\label{prop:r-valued2}
For $h\in \Def(X)$,
\begin{eqnarray}
\label{eq:r-valued2}
    \int_X h\,\dchifloor
    &=&
    \lim_{\epsilon\to 0^+}\frac{1}{\epsilon}\int_{\real}
    s\,\chi\{s\leq h<s+\epsilon\}\,ds
\\
    \int_X h\,\dchiceil
    &=&
    \lim_{\epsilon\to 0^+}\frac{1}{\epsilon}\int_{\real}
    s\,\chi\{s< h\leq s+\epsilon\}\,ds .
\end{eqnarray}
\end{proposition}
{\em Proof:} For $h$ affine on an open $k$-simplex $\sigma$, and
$0<\epsilon$ sufficiently small, $\int_{\real} s\,\chi\{s\leq
h<s+\epsilon\}\,ds =
\epsilon\,(-1)^k\left(-\frac{\epsilon}{2}+\inf_{\sigma} h\right)$ and
$\int_{\real} s\,\chi\{s< h \leq s+\epsilon\}\,ds =
\epsilon\,(-1)^k\left(-\frac{\epsilon}{2}+\sup_{\sigma} h\right)$. \qed

\subsection{Morse theory}
\label{sec:morse}

One important indication that the definition of $\int \dchifloor$ is
correct for our purposes is the natural relation to Morse theory:
the integrals with respect to $\dchifloor$ and $\dchiceil$ are Morse
index weighted sums of critical values of the integrand. This is a
localization result, reducing from an integral over all of $X$ to an
integral over an often discrete set of critical points.

Recall that a $C^2$ function $h:M\to\real$ on a smooth manifold $M$
is \style{Morse} if all critical points of $h$ are nondegenerate, in
the sense of having a nondegenerate Hessian matrix of second partial
derivatives. Denote by $\crit(h)$ the set of critical points of $h$.
For each $p\in\crit(h)$, the \style{Morse index} of $p$, $\mu(p)$,
is defined as the number of negative eigenvalues of the Hessian at
$p$, or, equivalently, the dimension of the unstable manifold $W^u(p)$ of the
vector field $-\nabla h$ at $p$.

Stratified Morse theory \cite{GM} is a powerful generalization to
triangulable spaces, including definable sets with respect to an
o-minimal structure \cite{BK,Schurmann}. We may interpret
$\dchifloor$ and $\dchiceil$ in terms of the  weighted stratified Morse
index of the graph of the integrand.

\begin{definition}
For $X\subset\real^n$ definable and $h\in\Def(X)$, define the
co-index of $h$, $\Index^*h$ to be the stratified Morse index of the
graph of $h$, $\Gamma_h\subset X\times\real$, with respect to the
projection $\pi:X\times\real\to\real$:
\begin{equation}
\label{eq:index}
    (\Index^*h)(x) = \lim_{\epsilon'\ll\epsilon\to 0^+}
        \chi\left(
        \overline{B_\epsilon(x)}
        \cap
        \{h<h(x)+\epsilon'\}
        \right) ,
\end{equation}
where $\overline{B_\epsilon(x)}$ is the closed ball of radius
$\epsilon$ about $x\in X$. The index $\Index_*$ is the stratified Morse
index with respect to height function $-\pi$: i.e., $\Index_*h =
\Index^*(-h)$ or
\begin{equation}
\label{eq:index2}
    (\Index_*h)(x) = \lim_{\epsilon'\ll\epsilon\to 0^+}
        \chi\left(
        \overline{B_\epsilon(x)}
        \cap
        \{h>h(x)-\epsilon'\}
        \right) .
\end{equation}
\end{definition}

Note that $\Index_*,\Index^*:\Def(X)\to CF(\overline{X})$, and the
restriction of these operators to $CF(X)$ is the identity (every point
of a constructible function is a critical point). The two types of
integration on $\Def(X)$ correspond to the Morse indices of the graph
with respect to the two orientations of the graph axis --- the
projections $\pi$ and $-\pi$.

\begin{theorem}
\label{thm:stratified}
For any continuous $h\in\Def(X)$,
\begin{equation}
\label{eq:morseindex}
        \int_X h\,\dchifloor = \int_{\overline{X}}h\Index^*h\,d\chi
         \quad  ;  \quad
        \int_X h\,\dchiceil   = \int_{\overline{X}}h\Index_*h\,d\chi
.
\end{equation}
\end{theorem}
{\em Proof:} On an open $k$-simplex $\sigma\subset X\subset\real^n$
in an affine triangulation of $h$, the co-index $\Index^*h$ equals
$(-1)^{\dim(\sigma)}$ times the characteristic function of the closed
face of $\sigma$ determined by $\inf_\sigma h$. Since $h$ is
continuous, $\int_{\overline{\sigma}}h\Index^*h\,d\chi =
(-1)^{\dim(\sigma)}\inf_\sigma h$. Lemma \ref{lem:r-simplex} and
additivity complete the proof; the analogous proof holds for
$\Index_*$ and $\dchiceil$. \qed

%

\begin{corollary}
\label{cor:morse}
If $h$ is a Morse function on a closed $n$-manifold $M$, then:
\begin{align}
\label{eq:morse}
    \int_M h\,\dchifloor
    =
    \sum_{p\in\crit(h)}(-1)^{n-\mu(p)}h(p)
    ;
\\
\label{eq:morse2}
    \int_M h\,\dchiceil
    =
    \sum_{p\in\crit(h)}(-1)^{\mu(p)}h(p)
    .
\end{align}
\end{corollary}
{\em Proof:} For $p\in\crit(h)$ a nondegenerate critical point on an
$n$-manifold, $\Index^*h(p)=(-1)^{n-\mu(p)}$ and
$\Index_*h(p)=(-1)^{\mu(p)}$. \qed

From this, one sees clearly that the relationship between $\dchifloor$
and $\dchiceil$ is regulated by Poincar\'e duality. For example, on
continuous definable integrands over an $n$-dimensional manifold
$M$,
\begin{equation}
\label{eq:morse3}
    \int_M h\,\dchiceil
    =
    (-1)^n\int_M h\,\dchifloor
    .
\end{equation}

The generalization from continuous to general definable integrands is
simple, but requires weighting $\Index^*h$ by $h$ directly. To
compute $\int_X h\dchifloor$, one integrates the weighted co-index
\begin{equation}
        \lim_{\epsilon'\ll\epsilon\to 0^+}
        h(x+\epsilon')\chi\left(
        \overline{B_\epsilon(x)}
        \cap
        \{h<h(x)+\epsilon'\}
        \right)
\end{equation}
with respect to $d\chi$.

Corollary \ref{cor:morse} can also be proved directly using classical
handle-addition techniques or in terms of the Morse complex, using
the fact that the restriction of the integrand to each unstable
manifold of each critical point is unimodal with a unique maximum at
the critical point. It is also possible to express the stratified Morse
index --- and thus the integral here considered --- in terms of
integration against a characteristic cycle, \cf \cite{GM,Schurmann}.

One final means of illustrating Corollary \ref{cor:morse} is to use a
deformation argument. Let $h$ be smooth on $X$ and $\phi_t$ be the
flow of $-\nabla h$. Then the integral is invariant under the action of
$\phi_t$ on $h$; yet the limiting function
$h_\infty=\lim_{t\to\infty}h\circ\phi_t$ is constant on stable
manifolds of $-\nabla h$ with values equal to the critical values of
$h$. We have not shown that the limiting function is constructible
(this depends on the existence of definable invariant manifolds --- we
are unaware of relevant results in the literature) and thus do not rely
on this method for proof but rather illumination.

\section{The integral operator}

We consider properties of the integral operator(s) on $\Def(X)$.

\subsection{Linearity}

One is tempted to apply all the standard constructions of sheaf theory
(as in \cite{Schapira,Schapira:tom}) to $\int_X:\Def(X)\to\real$.
However, our formulation of the integral on $\Def(X)$ has a glaring
disadvantage.

\begin{lemma}
\label{lem:nonlinear} $\int_X:\Def(X)\to\real$ (via $\dchifloor$ or
$\dchiceil$) is not a homomorphism for $\dim X>0$.
\end{lemma}
{\em Proof:}
\[
    1 =
    \int_{[0,1]} 1\,\dchifloor
    \neq
    \int_{[0,1]} x\,\dchifloor + \int_{[0,1]} (1-x)\,\dchifloor
    = 1+1 = 2.
\]
\qed

This loss of functoriality can be seen as due to the fact that $\lfloor
f+g \rfloor$ agrees with $\lfloor f \rfloor + \lfloor g \rfloor$ only up to
a set of Lebesgue measure zero, not $\chi$-measure zero. The
nonlinear nature of the integral is also clear from Eqn.
[\ref{eq:morseindex}], as Morse data is non-additive.

\subsection{The Fubini Theorem}

In one sense, the change of variables formula trivializes (Lemma
\ref{lem:cov}). The more general change of variables formula
encapsulated in the Fubini theorem does not, however, hold for
non-constructible integrands.

\begin{corollary}
\label{cor:fubinifail}
The Fubini theorem fails on $\Def(X)$ in general.
\end{corollary}

{\em Proof:} Let $F:X=Y\coprod Y\to Y$ be the
projection map with fibers $\{p\}\coprod\{p\}$. Any
$h\in\Def(X)$ is expressible as $f\coprod g$ for $f,g\in\Def(Y)$. The
Fubini theorem applied to $F$ is equivalent to the statement
\[
    \int_Y f + \int_Y g = \int_X h = \int_Y F_*h = \int_Y f+g
\]
(where the integration is with respect to $\dchifloor$ or $\dchiceil$ as
desired). Lemma \ref{lem:nonlinear} completes the proof. \qed

Fubini holds when the map respects fibers.

\begin{theorem}
\label{thm:r-fubini}
For $h\in\Def(X)$, let $F:X\to Y$ be definable and $h$-preserving
($h$ is constant on fibers of $F$). Then $\int_Y F_*h\dchifloor =
\int_X h\dchifloor$, and $\int_Y F_*h\dchiceil = \int_X h\dchiceil$.
\end{theorem}
{\em Proof:} An application of the o-minimal Hardt theorem
\cite{vdD} implies that $Y$ has a partition into definable sets
$Y_\alpha$ such that $F\inv(Y_\alpha)$ is definably homeomorphic to
$U_\alpha\times Y_\alpha$ for $U_\alpha$ definable, and that
$F:U\times Y_\alpha\to Y_\alpha$ acts via projection. Since $h$ is
constant on fibers of $F$, one computes
\[
    \int_{Y_\alpha}F_*h\dchifloor
    =
    \int_{Y_\alpha}h\,\chi(U_\alpha)\dchifloor
    =
    \int_{U_\alpha\times Y_\alpha}h\dchifloor .
\]
The same holds for $\int\dchiceil$. \qed

\begin{corollary}
\label{cor:r-fubini}
For $h\in\Def(X)$, $\int_Xh = \int_{\real}h_*h$. In other words,
\begin{equation}
    \int_Xh\,\dchifloor
    =
    \int_\real s\,\chi\{h=s\}\dchifloor ,
\end{equation}
and likewise for $\dchiceil$.
\end{corollary}

\subsection{Continuity}

Though the integral operator is not linear on $\Def(X)$, it does
retain some nice properties. All properties below stated for
$\int\,\dchifloor$ hold for $\int\,\dchiceil$ via duality.

\begin{lemma}
\label{lem:homog}
The integral $\int\,\dchifloor:\Def(X)\to\real$ is positively
homogeneous.
\end{lemma}
{\em Proof:} For $f\in \Def(X)$ and $\lambda\in\real^+$, the change
of variables variables $s\mapsto \lambda s$ in Eqn. [\ref{eq:r-valued}]
gives $\int\lambda f\,\dchifloor = \lambda\int f\,\dchifloor$. \qed

%
%

Integration is not continuous on $\Def(X)$ with respect to the $C^0$
topology. An arbitrarily large change in $\int h\dchifloor$ may be
effected by small changes to $h$ on a (large) finite point set.
\def\strat{{\mathcal S}}
\def\kk{{\mathbf k}}
In some situations the ``complexity'' of the definable functions can
be controlled in a way sufficient to ensure continuity.

One example arises in the semialgebraic category. Fix a (finite)
semialgebraic stratification $\strat$ of a compact definable $X$, and
consider definable semialgebraic functions {\em algebraic} with
respect to this stratification (that is such that the restriction of the
function to any stratum $S\in\strat$ is a polynomial $P_S$). The
resulting linear space (filtered by the subspaces of polynomials of
bounded degree) can be equipped with the structure of a Banach
space, by completing the family of semi-norms
$\|P\|_{S,\kk}=\max_{S\in \strat}||P_S||_{C^n}$, where $n=\dim X$.
Then $\int_X\cdot\dchifloor$ becomes a continuous (non-linear)
functional on this Banach space. The proof results, essentially, from
the B\'ezout theorem (mimicking Thom-Milnor theory): the total
number of critical points graph of a polynomial of degree $D$ on a
fixed semi-algebraic set is bounded by $O(D^n)$. The generalization
to increasing (refined) stratifications is straightforward.

Integration itself defines a natural topology for $\Def(X)$ on which
integration is continuous. Define the $L^1$ $\epsilon$-neighborhood
of $h\in\Def(X)$ as the intersection of the $C^0$
$\epsilon$-neighborhood (definable functions with $\epsilon$-close
graphs) with those functions $g\in\Def(X)$ satisfying $\abs{\int_X
f-g\,\dchifloor} < \epsilon$. This provides a basis for an $L^1$
topology on $\Def(X)$. As a consequence of Lemma
\ref{lem:floortoceiling}, the definition is independent of the use of
$\dchifloor$ or $\dchiceil$.

The interested reader may speculate on other function space
topologies on $\Def(X)$.

\subsection{Duality and links}
\label{sec:duality}
There is an integral transform on $CF(X)$ that is the analogue of
Poincar\'e-Verdier duality \cite{Schurmann}. It extends seamlessly
to integrals on $\Def(X)$ by means of the following definition.

\begin{definition}
\label{def:dual}
The \style{duality operator} is the integral transform $\Dual:CF(X)\to
CF(X)$ given by
\begin{equation}
\label{eq:dual}
    \Dual h(x)
    =
    \lim_{\epsilon\to 0^+}\int_Xh\one_{B_\epsilon(x)}d\chi ,
\end{equation}
where $B_\epsilon$ is an open metric ball of radius $\epsilon$.
\end{definition}

We extend the definition to $\Dual:\Def(X)\to\Def(X)$ by integrating
with respect to $\dchifloor$ or $\dchiceil$, interchangeably, via:

\begin{lemma}
\label{lem:dualindependent}
$\Dual h$ is well-defined on $\Def(X)$ and independent of whether
the integration in (\ref{eq:dual}) is with respect to $\dchifloor$
or $\dchiceil$.
\end{lemma}
{\em Proof:} The limit is always well-defined thanks to the Conic
Theorem in o-minimal geometry \cite{vdD}. To show that it is
independent of the upper- or lower-semicontinuous approximation,
take $\epsilon>0$ sufficiently small. Note that by triangulation, $h$
can be assumed to be piecewise-affine on open simplices. Pick a point
$x$ in the support of $h$ and let $\{\sigma_i\}$ be the set of open
simplices whose closures contain $x$. Then for each $i$, the limit
$h_i(x):=\lim_{y\to x}h(y)$ for $y\in\sigma_i$ exists. Then, applying
Eqn. [\ref{eq:dual}], one computes
\begin{equation}
\label{eq:dualsimplices}
    \Dual h(x)
    =
        \sum_i (-1)^{\dim\sigma_i}h_i(x) ,
\end{equation}
independent of the measure $\dchifloor$ or $\dchiceil$. \qed

For a continuous definable function $h$ on a manifold $M$, $\Dual h =
(-1)^{\dim M}h$, as one can verify by combining Eqns.
[\ref{eq:r-valued}] and [\ref{eq:dual}]. This is commensurate with the
result of Schapira \cite{Schapira} that $\Dual$ is an involution on
$CF(X)$.

\begin{theorem}
Duality is involutive on $\Def(X)$: $\Dual\circ\Dual h = h$.
\end{theorem}
{\em Proof:} Given $h$, fix a triangulation on which $h$ is affine on
open simplices. From Eqn. [\ref{eq:dualsimplices}], we see that the
dual of $h$ at $x$ is completely determined by the trivialization of
$h$ at $x$. Let $L_x{h}$ be the constructible function on
$B_\epsilon(x)$ which takes on the value $h_i(x)$ on strata
$\sigma_i\cap B_\epsilon(x)$. (Though this is not necessarily an
integer-valued function, its range is discrete and therefore it is
constructible.) As $L_xh$ is close to $h$ in $B_\epsilon(x)$ (this
follows from the continuity of $h$ on each of the strata), $\Dual h$ is
close to $\Dual L_xh$ in $B_\epsilon(x)$: indeed, the total Betti
number of intersections of strata with any ball $B_\epsilon(y)$ is
bounded, and Euler integral of a function small in absolute value is
small as well. Hence the definable function $\Dual^2 h$ is close to
the constructible function $\Dual^L_xh$ with $\epsilon$ small. As
$\Dual^2L_x{h}(x)=L_x{h}(x)=h(x)$, the result follows. \qed

One can define related integral transforms. For example, the
\style{link} of $h\in CF(X)$ is defined as
\begin{equation}
\label{eq:link}
    \Lambda h(x)
    =
    \lim_{\epsilon\to 0^+}\int_Xh\one_{\partial B_\epsilon(x)}d\chi.
\end{equation}
The link of a continuous function on an $n$-manifold $M$ is
multiplication by $1+(-1)^n$, as a simple computation shows. In
general, $\Lambda = \id - \Dual$, where $\id$ is the identity
operator.

\section{Integral transforms}

Integration with respect to Euler characteristic over $CF(X)$ has a
well-defined and well-studied class of integral transforms, expressed
beautifully in Schapira's work on inversion formulae for the
generalized Radon transform in $d\chi$ \cite{Schapira:tom}. Integral
transforms with respect to $\dchifloor$ and $\dchiceil$ are similarly
appealing, with applications to signal processing as a primary
motivation. Examples of interesting definable kernels for integral
transforms over Euclidean $\real^n$ include $\|x-y\|$, $\langle x,y
\rangle$, and $g(x-y)$ for some $g$. These evoke Bessel (Hankel) transforms,
Fourier transforms, and convolution with $g$ respectively. The choice
between $\dchifloor$ and $\dchiceil$ makes a difference, of course,
but can be amalgamated. Example: for fixed kernel $K$, one
can consider the mixed integral transform $h\mapsto \int_X hK \dchifloor -
\int_X hK\dchiceil$. In the case of  $K(x,\xi)=\langle x,\xi \rangle$,
this transform takes $\one_A$ for $A$ compact and convex to the
`width' of $A$ projected to the $\xi$-axis.

\subsection{Convolution}
\label{sec:convolution}

On a vector space $V$ (or Lie group, more generally), a convolution
operator with respect to Euler characteristic is straightforward. Given
$f,g\in CF(V)$, one defines
\begin{equation}
\label{eq:conv}
    (f*g)(x) = \int_V f(t)g(x-t)\,d\chi .
\end{equation}

Convolution behaves as expected in $CF(V)$. By reversing the order of
integration, one has immediately that $\int_V f*g \,d\chi = \int_V
f\,d\chi \, \int_V g\,d\chi$. There is a close relationship between
convolution and the Minkowski sum, as observed in, \eg,
\cite{Groemer}: for $A$ and $B$ convex and closed
$\one_A*\one_B=\one_{A+B}$, cf. \cite{Viro,Schapira}. Convolution is
a commutative, associative operator providing $CF(V)$ with the
structure of an (interesting \cite{Brocker}) algebra.

Convolution is well-defined on $\Def(V)$ by integrating with respect
to $\dchifloor$ or $\dchiceil$. However, the product formula for $\int
f*g$ fails in general, since one relies on the Fubini theorem to prove it
in $CF(V)$.

\subsection{Linearity}

The nonlinearily of the integration operator prevents most
straightforward applications of inversion formulae \`a la Schapira. Fix
a kernel $K\in\Def(X\times Y)$ and consider the integral transform
$\Transform_K:\Def(X)\to\Def(Y)$ of the form $(\Transform_K
h)(y)=\int_Xh(x)K(x,y)\dchifloor(x)$. In general, this operator is
non-linear, via Lemma \ref{lem:nonlinear}. However, some vestige of
(positive) linearity survives within $CF$.

\begin{lemma}
\label{lem:CF^+}
The integral transform $\Transform_K$ is positive-linear over
$CF^+(X)=\Def(X,\nats)$.
\end{lemma}

{\em Proof:} Any $h\in CF^+(X)$ is of the form $h=\sum_k
a_k\one_{U_k}$ for $a_k\in\nats$ and $U_k\in\Def(X)$. For
$h=\one_A$, $\Transform_K h = \int_A K \dchifloor$. Additivity of the
integral in $\dchifloor$ (via Eqn. [\ref{eq:measure}]) combined with
Lemma \ref{lem:homog} completes the proof. \qed

This implies in particular that when one convolves a function $h\in
CF^+(\real^n)$ with a smoothing kernel (e.g., a Gaussian) as a means
of filtering noise or taking an average of neighboring data points, that
convolution may be analyzed one step at a time (decomposing $h$).

Integral transforms are not linear over all of $CF(X)$, since $\int
-h\dchifloor \neq -\int h\dchifloor$. However, integral transforms
which combine $\dchifloor$ and $\dchiceil$ compensate for this
behavior. Define the measure $[d\chi]$ to be the average of $\dchifloor$
and $\dchiceil$:
\begin{equation}
    \int_X h[d\chi] = \frac{1}{2}\left(\int_X h\dchifloor + \int_X h\dchiceil\right) .
\end{equation}

\begin{theorem}
\label{thm:lineartransform}
Any integral transform of the form
\begin{equation}
    (\Transform_K h)(y)
    =  \int_X h(x)K(x,y) [d\chi](x)
\end{equation}
is a linear operator $CF(X)\to\Def(Y)$.
\end{theorem}
{\em Proof:} From Lemma \ref{lem:CF^+}, $\Transform$ is
positive-linear over $CF^+(X)$. Full linearity follows from the
observation that $\int_{X} -h [d\chi]  = -\int_{X} h [d\chi]$, which
follows from Lemma \ref{lem:r-simplex} by triangulating $h$. \qed

As a simple example, consider the transform with kernel
$K(x,\xi)=\langle x,\xi\rangle$. The transform of $\one_A$ with
respect to $[d\chi]$ for $A$ compact and convex equals a `centroid'
of $A$ along the $\xi$-axis: the average of the maximal and minimal
values of $\xi$ on $\partial A$. Note how the dependence on critical
values of the integrand on $\partial A$ reflects the Morse-theoretic
interpretation of the integral in this case.

Integration with respect to $[d\chi]$ seems suitable only for
integral transforms over $CF$. On a continuous integrand, the
integral with respect to $[d\chi]$ either returns zero (cf. the
integral of Rota \cite{Rota}) or else the integral with respect to
$\dchifloor$, depending on the parity of the $\dim X$, via Eqn.
[\ref{eq:morse3}].

\section{Applications of definable Euler integration}
\label{sec:app}

The Euler calculus on $CF$ is quite useful; the extension to $\Def$
deepens this utility and opens new potential applications, of which
we highlight a few.

\subsection{Sensor networks}
\label{sec:app:sn}

An application of Euler integration over $CF(X)$ to sensor networks
problems was initiated in \cite{BG:enum}. Consider a space $X$ whose
points represent target-counting sensors that scan a workspace $W$.
Target detection is encoded in a \style{sensing relation}
$\Support\subset W\times X$ where $(w,x)\in\Support$ iff a target at
$w$ is detected by a sensor at $x$. Assume that sensors count the
number of sensed targets, but do not locate or identify the targets.
The sensor network therefore induces a \style{target counting
function} $h:X\to\nats$ of the form $h=\sum_\alpha\one_{U_\alpha}$,
where $U_\alpha$ is the \style {target support} --- the set of
sensors which detect target $\alpha$. Euler integration allows for
simple enumeration:

\begin{theorem}[\cite{BG:enum}]
Assume $h\in CF(X)$ and $\chi(U_\alpha)=N\neq 0$ for all $\alpha$.
Then the number of targets in $W$ is precisely $\frac{1}{N}\int_X
h\,d\chi$.
\end{theorem}

Since the target count is presented as an integral, it is possible
to accurately estimate the answer when the integrand $h$ is known
not on all of $X$ (a continuum of sensors being an idealization) but
rather on a sufficiently dense grid of sample points (physical
sensors in a network).

The $\real$-valued theory aids in establishing expected values of
target counts in the presence of confidence measures on sensor
readings. Let $\Nodes=\{x_i\}$ denote a discrete set of sensor nodes
in $\real^n$, and assume each sensor returns a target count
$h(x_i)\in\nats$ and a fluctuation measure $c(x_i)\in[0,1]$ obtained,
say, by stability of the reading over a time average. View $h$ as a
sampling over $\Nodes$ of the true target count
$f=\sum_\alpha\one_{U_\alpha}$. Assume that nodes with fluctuation
reading $0$ have perfect information ($h=f$ at $x_i$) and that $c$
correlates with error $\abs{f-h}$. Assume that sensor nodes $\Nodes$
are the vertex set of a triangulation $\Triangulation$.

The integral of an extension of $f$ over a triangulation gives a
terrible approximation to $\int h\,d\chi$: an error of $\pm 1$ on $K$
nodes can cause a change in the integral of order $K$. More
specifically, if $h=f+e$, where $e:\Nodes\to\{-1,0,1\}$ is an error
function that is nonzero on a sparse subset $\Nodes'\subset\Nodes$,
then, for certain infelicitous choices of $\Nodes'$, $\abs{\int h - \int
f} = \abs{\Nodes'}$.

A $\real$-valued relaxation can mitigate errors by using fluctuation
$c$ as a weight. Let $N(i)$ be a collection of neighboring nodes to
$x_i$, where neighborhood can be defined via distance (if available)
or edge-distance (in an ad hoc network or triangulation). Define
$\tilde{h}$ to be the result of averaging the value at $x_i\in\Nodes$
over $N(i)$, with $c$ as a weight. Specifically,
\begin{equation}
    \tilde{h}(x_i) = \frac
    {\sum_{y\in N(i)}c(y)h(y)}
    {\sum_{y\in N(i)}c(y)}.
\end{equation}
This nearest-neighbor convolution damps out local variations. The
resulting integral with respect to $\dchifloor$ will tend to mitigate
localized errors, thanks to the Morse-theoretic formula: see Figs. 1-2.
More numerical investigation is warranted.

\begin{figure}
\begin{center}
\includegraphics[angle=0,width=3.0in]{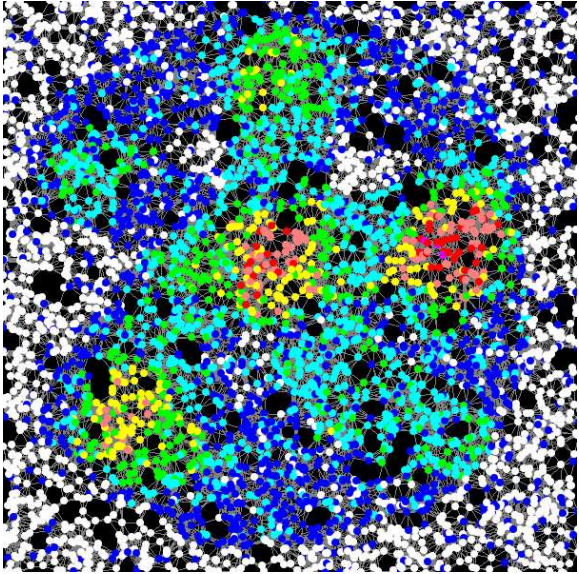}
\caption{Shown is a sampling of an upper semicontinuous
constructible integrand $h:\real^2\to\nats$, sampled over a random
network with holes. Values of the integrand are color-coded by
height, white being zero. One third of the sample points sample $h$
with an error of $\pm 1$, uniformly distributed. The Euler integral
of this sampling with respect to $d\chi$ is $64$, a poor
approximation to the true $\int h\, d\chi = 9$.}
\end{center}
\end{figure}
%
\begin{figure}
\begin{center}
\includegraphics[angle=0,width=3.0in]{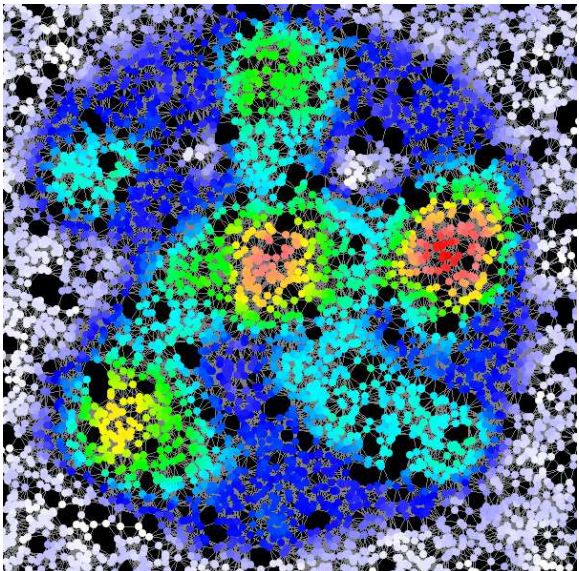}
\caption{The above sampled integrand, averaged over nearest
neighbors. The integral of this $\real$-valued integrand with
respect to $\dchifloor$ is approximately $9.52$: a reasonable
approximation to the true $\int h\, d\chi = 9$.}
\end{center}
\end{figure}

Such averaging leads to non integer-valued integrands. By using
integration with respect to $\dchifloor$ or $\dchiceil$ for
upper/lower semi-continuous integrands associated to an averaged
signal $\tilde{h}$, one obtains an expected value of $\int
h\,d\chi$. This can be particularly illuminating when a network has
incomplete information, \eg, a hole. Holes in a network can be
modeled by setting the confidence measure $c$ to zero and averaging.

\subsection{Statistics and mode counting}
\label{sec:app:mc}

The previous applications lend themselves to more general
statistical ends. Consider a smooth distribution
$f:\real^n\to[0,\infty)$ of compact support and bounded variation.
The problem of mode-counting --- of decomposing $f$ into a convex
combination of unimodal summands --- is of interest to
statisticians, even in the univariate case \cite{E,B}. Gaussian
summands are commonly used, though not exclusively \cite{Ka,Ke}. The
techniques of this paper are relevant (see \cite{RL} for Morse
structures associated to mixtures of multivariate normal
distributions).

One way to interpret mode-counting is as a topological deconvolution
problem. Assume that $f$ is of the form $f=u*h$ for  $h\in
CF^+(\real^n)$ and $u$ a unimodal distribution supported on a
neighborhood of $0\in\real^n$, where \style{unimodal} means that all
upper excursion sets $\{u>c\}$ are contractible for $c>0$. One
assumes that $h$ is of the form $h=\sum_{i=1}^n\one_{U_i}$ where
the supports $U_i$ are sufficiently simple so that $u*\one_{U_i}$ has
no self-interference (the support of $u*\one_{U_i}$ strongly
deformation retracts to $U_i$).

The mode-counting problem is equivalent to computing $\int h\,d\chi$
given $f$ and some information about $u$, say, its height $\max u =
\int u\dchifloor$. If the convolution formula for $CF$ held over
$\Def$, then the number of modes would be $\int f\dchifloor$ divided
by $u_{\max}$. However, the non-linearity of the integral with
respect to $\dchifloor$ precludes this solution. Indeed, this
nonlinearity mirrors the interaction of unimodal summands in a
distribution. Just as two modes can interfere, creating a local
maximum when an increasing and a decreasing portion of the modes
are summed, the Euler integral over $\Def$ loses linearity when
increasing and decreasing integrands are combined.

The development of good algorithms for integral transforms over
$\dchifloor$ or even $[d\chi]$ will be useful not only for
mode-counting, but also for explicit mode decomposition, since the
extraction of the modes themselves from the unknown factor $h$
involves deconvolving $f=u*h$.

%

\subsection{Numerical integration}

Though integration with respect to Euler characteristic has a
lengthy history, there appears to be no treatment of numerical
integration, even in the simpler setting of $CF^+(\real^n)$.
The central problem (in the constructible and definable categories) is
how to estimate $\int_X h$ given the values of $h$ on a discrete
subset of $X$. As in the case of numerical integration for Riemann
integrals, one typically assumes something about the features of $h$
and/or the density and extent of the sampling set. In \cite{BG:enum},
the present authors give a formula for estimating $\int h\,d\chi$ given
a discrete sampling of $h\in CF(\real^2)$ which correctly samples
connectivity data of excursion sets. This formula generalizes to the
definable category:

\begin{proposition}
\label{prop:betti0}
For $h\in\Def(\real^2)$ continuous, $\int h\,\dchifloor =$
\begin{equation}
    \int_{s=0}^\infty
        \beta_0\{h\geq s\}
        +
        \beta_0\{h\geq-s\}
        -
        \beta_0\{h<s\}
        -
        \beta_0\{h<-s\}
        \,ds,
\end{equation}
where $\beta_0(\cdot)=\dim H_0(\cdot;\real)$, the zeroth Betti
number.
\end{proposition}
{\em Proof:} Apply the homological definition of $\chi$ to Eqn.
[\ref{eq:r-valued}]; then, use Alexander duality in the plane to reduce
all terms to $\beta_0$ quantities. \qed

The value of Proposition \ref{prop:betti0} is that it allows for
computation based on $\beta_0$ quantities. Such connectivity data
are easily obtained from a discrete sampling via clustering. We have
implemented this formula in software (see Fig. \ref{fig:ave4}).
However, for more general integration domains than $\real^2$,
duality formulae are less helpful. One general result on refinement
follows from continuity of the integral operator.

\begin{theorem}
\label{thm:cont}
For $h\in\Def(X)$ continuous, let $h_{PL}$ be the piecewise-linear
function obtained from sampling $h$ on the vertex set of a
triangulation $\Triangulation$ of $X$. As the sampling and
triangulation are refined,
\begin{equation}
    \lim_{\abs{\Triangulation}\to 0^+}\int_Xh_{PL}\dchifloor
    =
    \int_Xh\dchifloor .
\end{equation}
\end{theorem}

This result relies crucially on continuity and does not apply to
$CF(X)$. A more desirable result would be a measure of how far a
given sampling is from the true integral. This seems challenging. We
note that the Morse-theoretic formulae
[\ref{eq:morse}]-[\ref{eq:morse2}] allow one to reduce the domain of
an integral to a (typically finite) set of critical points. This `focusing'
property of integration over $\Def$ should be a starting point for good
numerical algorithms, especially for integral transforms.

\begin{acknowledgments}
This work supported by DARPA \# HR0011-07-1-0002 and ONR
N000140810668.
\end{acknowledgments}

\end{article}









\end{document}